\documentclass[11pt,twoside,reqno]{amsart}
\usepackage{import}
\usepackage{newclude}
\usepackage{anysize}
\usepackage{amsmath}
\usepackage{amsthm}
\usepackage{amsmath,amscd}
\usepackage[utf8]{inputenc}
\usepackage{amssymb}
\usepackage{stmaryrd}
\usepackage{wasysym}
\usepackage{mathrsfs}
\usepackage{hyperref}
\usepackage{cleveref}
\usepackage{enumitem}
\usepackage{dsfont}
\usepackage{soul}
\hypersetup{colorlinks=true,linkcolor=black,citecolor=black}
\usepackage{color}
\usepackage[all]{xy}
\usepackage{float}
\usepackage{bm}
\usepackage{mathtools}
\usepackage{thmtools}
\usepackage{setspace}
\usepackage{comment}
\usepackage{tikz}
\usepackage{tikz-cd}
\usepackage{mathdots}
\usepackage{graphicx}
\usepackage[myheadings]{fullpage}

\numberwithin{equation}{section}
\newcommand\mtop{.95in}
\newcommand\mbottom{.95in}
\newcommand\mleft{1in}
\newcommand\mright{1in}
\usepackage[top = \mtop, bottom = \mbottom, left = \mleft, right=\mright]{geometry}

\DeclareMathOperator{\Mat}{Mat}

\newtheorem{thm}{Theorem}[section]
\newtheorem{example}[thm]{Example}
\newtheorem{prop}[thm]{Proposition}

\newtheorem{lemma}[thm]{Lemma}

\newtheorem{cor}[thm]{Corollary}
\theoremstyle{definition}
\newtheorem{defi}[thm]{Definition}

\newtheorem{rmk}[thm]{Remark}

\newcommand\reallywidehat[1]{%
\savestack{\tmpbox}{\stretchto{%
  \scaleto{%
    \scalerel*[\widthof{\ensuremath{#1}}]{\kern-.6pt\bigwedge\kern-.6pt}%
    {\rule[-\textheight/2]{1ex}{\textheight}}
  }{\textheight}%
}{0.5ex}}%
\stackon[1pt]{#1}{\tmpbox}%
}
\DeclareSymbolFont{bbold}{U}{bbold}{m}{n}
\DeclareSymbolFontAlphabet{\mathbbold}{bbold}

\makeatletter
\def\@tocline#1#2#3#4#5#6#7{\relax
  \ifnum #1>\c@tocdepth 
  \else
    \par \addpenalty\@secpenalty\addvspace{#2}%
    \begingroup \hyphenpenalty\@M
    \@ifempty{#4}{%
      \@tempdima\csname r@tocindent\number#1\endcsname\relax
    }{%
      \@tempdima#4\relax
    }%
    \parindent\z@ \leftskip#3\relax \advance\leftskip\@tempdima\relax
    \rightskip\@pnumwidth plus4em \parfillskip-\@pnumwidth
    #5\leavevmode\hskip-\@tempdima
      \ifcase #1
       \or\or \hskip 1em \or \hskip 2em \else \hskip 3em \fi%
      #6\nobreak\relax
    \hfill\hbox to\@pnumwidth{\@tocpagenum{#7}}\par
    \nobreak
    \endgroup
  \fi}
\makeatother

\newcommand{\R}{\mathbb{R}}
\newcommand{\Z}{\mathbb{Z}}
\newcommand{\Q}{\mathbb{Q}}
\newcommand{\N}{\mathbb{N}}

\newcommand{\C}{\mathbb{C}}

\newcommand{\E}{\mathbb{E}}

\renewcommand{\l}{\lambda}

\newcommand{\var}{\text{Var}}

\renewcommand{\L}{\Lambda}

\newcommand{\G}{\mathbb{G}}

\newcommand{\K}{K}

\DeclareMathOperator{\Hom}{Hom}

\DeclareMathOperator{\GL}{GL}
\DeclareMathOperator{\SL}{SL}

\DeclareMathOperator{\Sp}{Sp}
\DeclareMathOperator{\GSp}{GSp}

\DeclareMathOperator{\Cov}{Cov}

\DeclareMathOperator{\SN}{SN}
\DeclareMathOperator{\Cor}{Cor}
\DeclareMathOperator{\diag}{diag}

\DeclareMathOperator{\id}{id}

\title{Gaussian Universality of Products Over Split Reductive Groups and the Satake Isomorphism}
\author{Jiahe Shen}

\date{\today}
\begin{document}

\thanks{I thank my advisor Ivan Corwin for helpful advice on revisions and for providing me funding support with his NSF grant DMS-2246576 and Simons Investigator grant 929852; Roger Van Peski, for suggestions on the topics, valuable encouragement, and comments for polishing my rhetorical skills; Professor Chao Li, for providing materials around split reductive groups, Satake isomorphism and related topics; I also wish to thank Yifan Wu, Fan Zhou for additional helpful conversations.}

\maketitle

\begin{abstract} 
We establish that the singular numbers (arising from Cartan decomposition) and corners (emerging from Iwasawa decomposition) in split reductive groups over non-archimedean fields are fundamentally determined by Hall-Littlewood polynomials. Through applications of the Satake isomorphism, we extend Van Peski's results \cite[Theorem 1.3]{van2021limits} to encompass arbitrary root systems. Leveraging this theoretical foundation, we further develop Shen's work \cite[Theorem 1.1]{shen2024gaussian} to demonstrate that both singular numbers and corners of such products exhibit minimal separation. This characterization enables the derivation of asymptotic properties for singular numbers in matrix products, particularly establishing the strong law of large numbers and central limit theorem for these quantities. Our results provide a unified framework connecting algebraic decomposition structures with probabilistic limit theorems in non-archimedean settings.
\end{abstract}



\tableofcontents

\section{Introduction}
\subsection{Asymptotic results}\label{subsec: Asymtotic results}

Random matrix theory occupies a central position in modern probability and mathematical physics. The asymptotic behavior of random matrix products over non-archimedean fields (e.g., 
$p$-adic fields) introduces algebraic and probabilistic features due to the group structures and decomposition theorems (e.g., Cartan/Iwasawa decomposition, etc.) inherent to these settings. Investigations in this domain explore the applications of split reductive groups within probability theory and provide insights into related stochastic processes.

Consider classical groups like the special linear group $G=\SL_{n+1}(F)$ or the symplectic group $G=\Sp_{2n}(F)$ over a non-archimedean field $F$. In these examples, the Cartan decomposition $G=K\L^+K$ defines \emph{singular numbers}, regarded as the analogs of eigenvalues for non-archimedean matrices. On the other hand, the Iwasawa decomposition, as the generalization of Gram-Schmidt orthogonalization, defines \emph{corners}. These invariants characterize the geometric and probabilistic properties of matrix products. It is natural to ask about the evolution of singular numbers $\l(k)=\SN(A_1\cdots A_k)$ and their correlations with corners, which describes the dynamics based on non-archimedean random matrix products.

Recent work by Van Peski \cite[Theorem 1.1]{van2021limits} established Gaussian universality laws for singular numbers in $\GL_n(F)$, linking their statistics to Hall-Littlewood polynomials. Following this result, Shen \cite{shen2024gaussian} proves Gaussian universality, where the matrices no longer need to be distributed with respect to the Hall-Littlewood measure. The method can be summarized as follows: it involves comparing the singular values of matrix products with the sum of matrix corners and demonstrating that their differences become negligible in the asymptotic regime. Consequently, the Gaussian universality of matrix products follows from applying the strong law of large numbers and the central limit theorem to the sum of matrix corners. Utilizing this technique, the article also establishes asymptotic results for the $\GSp_{2n}(F)$ case. However, extending such results to general split reductive groups becomes challenging due to the intricate structure of root systems. This work addresses this gap by using representation-theoretic methods to analyze the probabilistic behavior of matrix products across different root systems in a unified way. We start with the following settings and keep them throughout this article:
\begin{enumerate}
\item (Root system and Weyl group)

Let $(V,\langle\cdot,\cdot\rangle)$ be a Euclidean space of dimension $n$;

$\Phi=\Pi^+\cup\Pi^-$ be a crystallographic root system that spans $V$, where $\Pi^+\subset\Phi$ is the subset of positive roots;

$\rho=\frac{1}{2}\sum_{\alpha\in\Pi^+}\alpha$ be the Weyl vector; 

$\alpha^\vee=\frac{2\alpha}{\langle\alpha,\alpha\rangle}$ be the coroot of $\alpha\in\Phi$, then we must have $\langle\alpha^\vee,\rho\rangle=1$;

$\Delta^+=\{\alpha_1,\ldots,\alpha_n\}$ be the corresponding set of positive simple roots;

$R=\sum_{1\le i\le n}\Z\alpha_i$ be the root lattice, which is the free abelian group generated by $\Delta^+$; 

$R^\vee=\sum_{1\le i\le n}\Z\alpha_i^\vee$ be the coweights, which is the free abelian group generated by coroots;

$R_0^\vee=\sum_{1\le i\le n}\N\alpha_i^\vee\subset R^\vee$ be the positive octant of $R^\vee$, where every coefficient has to be non-negative. This gives a partial ordering over $R^\vee$: for all $\l,\mu\in R^\vee$, we write $\l\ge\nu$ when $\l-\nu\in R_0^\vee$;

$R_+^\vee=\{\l\in R^\vee\mid\langle\l,\alpha_i\rangle\in\Z_{\ge 0},\forall 1\le i\le n\}\subset R_0^\vee$
be the dominant coweights; 

$W$ be the Weyl group generated by the elements $\{s_\alpha\in O(V)\mid\alpha\in\Phi\}$ where $s_\alpha:\l\mapsto\l-\frac{2\langle\l,\alpha\rangle}{\langle\alpha,\alpha\rangle}\alpha$ is the reflection that fix the hyperplane perpendicular to $\alpha$. For any $\mu\in R^\vee$, there exists a unique element $\l\in R_+^\vee$ such that $u=w\l$ for some $w\in W$. 

\item(Split reductive group)

Let $F$ be a non-archimedean local field whose residue field has finite order $q<\infty$, i.e., any finite algebraic
extensions of $\Q_p$ or $\mathbb{F}_p((t))$. We also denote $\mathfrak{o}\subset F$ as the ring of integers, $\mathfrak{p}$ as the maximal ideal, and $\pi$ as a generator of $\mathfrak{p}$. In particular, when $F=\Q_p$, we have $\mathfrak{o}=\Z_p,\mathfrak{p}=p\Z_p,q=|\Z_p/p\Z_p|=p$, and we may set $\pi=p$;

$G=G(F)$ be a split reductive group over $F$, with the Weyl group $W$ generated by the reflection of roots defined above;

$T$ be the maximal torus of $G$, which must be commutative of the form $\G_m^n$, the algebraic torus of dimension $n$. We have canonical isomorphism $R^\vee\cong\Hom(\G_m,T)$ that identifies $R^\vee$ as the coweight group of $T$.

$\L=\{\pi_\l\mid\l\in R^\vee\}$ be the subgroup of $T$, where $\pi_\l$ is the image of $\pi$ under $\l$, viewed as the homomorphism in $\Hom(\G_m,T)$;

$\L_+=\{\pi_\l\mid\l\in R_+^\vee\}\subset\L$ be the subset of $\L$, where the coweight is required to be dominant;

$B$ be the Borel subgroup of $G$, which must be solvable;

$N=B'$ be the unipotent radical of $B$, which leads to the semidirect product decomposition $B=T\rtimes N$. On the other hand, we have $T\cong B/N$;

$\K=G(\mathfrak{o})$ be the maximal compact subgroup of $G$, which gives the Iwasawa decomposition $G=BK= N\L K$;

$\mathfrak{g}=\mathfrak{t}\oplus\bigoplus_{\alpha\in\Phi}\mathfrak{g}_\alpha$ be corresponding weight space decomposition. Here $\mathfrak{t}$ corresponds to the maximal torus $T$, $\bigoplus_{\alpha\in\Pi^+}\mathfrak{g}_\alpha$ corresponds to the unipotent radical $N$, and $\mathfrak{t}\oplus\bigoplus_{\alpha\in\Pi^+}\mathfrak{g}_\alpha$ corresponds to the Borel subgrouop $B$.

\end{enumerate}

On the one hand, by the Cartan decomposition given in \cite[Proposition 4.4.3 (2)]{bruhat1972groupes}, for any $A\in G$, there exists $U,V\in K$ such that $A=U\pi_\l V$, where $\l\in R_+^\vee$ does not depend on our choice of $U,V$. We will denote $\SN(A)=\l$ as the \emph{singular numbers} of $A$. In particular, $\SN(A)=0$ if and only if $A\in K$. On the other hand, by the Iwasawa decomposition given in \cite[Proposition 4.4.3 (1)]{bruhat1972groupes}, there exists $U\in N, V\in K$ such that $A=U\pi_\nu V$, where $\nu\in R^\vee$ does not depend on our choice of $U,V$. We will denote $\Cor(A)=\nu$ as the \emph{corners} of $A$. The following example illustrates these settings and decompositions.

\begin{example}\label{ex: Cartan and Iwasawa for SL}
Consider the case $G=\SL_{n+1}(F)$, which is the type $A$ case. The Weyl group $W=S_{n+1}$ is the symmetric group that acts naturally on $V=\{(x_1,\ldots,x_{n+1})\subset\R^{n+1}\mid x_1+\cdots+x_{n+1}=0\}$. Let $\epsilon_i,1\le i\le n+1$ denote the standard basis of $\R^{n+1}$, then we have
$$\Pi^+=-\Pi^-=\{\epsilon_i-\epsilon_j\mid1\le i<j\le n+1\},\rho=\frac{1}{2}(n,n-2,\ldots,2-n,-n)$$
$$\Delta^+=\{\epsilon_i-\epsilon_{i+1}\mid1\le i\le n\},R=R^\vee=\{(\l_1,\ldots,\l_{n+1})\in\Z^{n+1}\mid\l_1+\cdots+\l_{n+1}=0\}$$
$$R^\vee_0=\{(\l_1,\ldots,\l_{n+1})\in\Z^{n+1}\mid\l_1+\cdots+\l_{n+1}=0,\quad\l_1+\cdots+\l_i\ge 0,\forall 1\le i\le n\}$$
$$R^\vee_+=\{(\l_1,\ldots,\l_{n+1})\in\Z^{n+1}\mid\l_1\ge\ldots\ge\l_{n+1},\l_1+\cdots+\l_{n+1}=0\}.$$
On the other hand, $T$ is the subset of $G$ consisting of diagonal matrices, $\L$ is the subset of $T$ where every entry on the diagonal is some power of $\pi$, $B$ (resp. $N$) is the subset of $G$ consisting of upper (resp. strict upper) triangular matrices, $K=\SL_{n+1}(\mathfrak{o})$. Let $A\in G=\SL_{n+1}(F)$. By Cartan decomposition, there exists $U,V\in\SL_{n+1}(\mathfrak{o})$ such that $A=U\diag_{(n+1)\times(n+1)}(\pi^{\l_1},\ldots,\pi^{\l_{n+1}})V$, where $(\l_1,\ldots,\l_{n+1})\in R^\vee_+$. Also, by the Iwasawa decomposition, there exists $U\in N,V\in\SL_{n+1}(\mathfrak{o})$ such that $A=U\diag_{(n+1)\times(n+1)}(\pi^{\nu_1},\ldots,\nu^{\nu_{n+1}})V$, where $(\nu_1,\ldots,\nu_{n+1})\in R^\vee$. Notice that the non-increasing relation $\nu_1\ge\ldots\ge\nu_{n+1}$ might not hold for the corners.
\end{example}


\begin{rmk}\label{rmk: the name corner}
Let us continue with the case $G=\SL_{n+1}(F)$. Here we shall explain why we name the coweight obtained from Iwasawa decomposition as ``corners''. For all $1\le i\le n+1$, let $A^{(i)}\in\Mat_{(n-i+2)\times(n+1)}(F)$ denote the corner of $A$ consisting of the last $n-i+2$ rows and all $n+1$ columns. Then there exists $U^{(i)}\in\GL_{(n-i+2)}(\mathfrak{o}),V^{(i)}\in\GL_{n+1}(\mathfrak{o})$ such that $A^{(i)}=U^{(i)}\diag_{(n-i+2)\times(n+1)}(\pi^{\l_1^{(i)}},\ldots,\pi^{\l_{n-i+2}^{(i)}})V^{(i)}$, where $\l_1^{(i)}\ge\ldots\ge\l_{n-i+2}^{(i)}$. Set $\SN(A^{(i)})=(\l_1^{(i)},\ldots,\l_{n-i+2}^{(i)})$, and $|\SN(A^{(i)})|=\l_1^{(i)}+\cdots+\l_{n-i+2}^{(i)}$ for all $1\le i\le n+1$.

Now, we claim that 
$$\Cor(A)=(|\SN(A^{(1)})|-|\SN(A^{(2)})|,\ldots,|\SN(A^{(n)})|-|\SN(A^{(n+1)})|,|\SN(A^{(n+1)})|)$$ 
and therefore $\Cor(A)$ provides information of the joint corners. Since 
$$
\Cor(A),\SN(A^{(1)}),\ldots,\SN(A^{(n+1)})$$ 
are all invariant under the right multiplication of $K$, there is no loss we assume $A\in N\L$ is upper triangular. In this case, $\Cor(A)=(\nu_1,\ldots,\nu_{n+1})$ comes from the entries on the diagonal. Then we must have $|\SN(A^{(i)})|=\nu_i+\ldots+\nu_{n+1}$, and thus our claim is verified.
\end{rmk}

Due to space limitations, we will not provide details for every root system, but one may also check \Cref{ex: Cartan and Iwasawa for Sp} for $G=\Sp_{2n}(F)$, which is the case of type $C$. 

Now, let us turn back to the general settings. Given random $A_1,A_2\ldots\in G$, we denote $\l(k)=\SN(A_1\cdots A_k)$ and $\nu(k)=\Cor(A_1\cdots A_k)$ for all $k\ge 1$. One may view $\l(1),\l(2),\ldots$ (resp. $\nu(1),\nu(2),\ldots$) as a discrete-time stochastic process on the lattice $R_+^\vee$ (resp. $R^\vee$). We prove the following theorem, which gives the limits and fluctuations of the sequence of singular numbers of products:

\begin{thm}\label{cor: i.i.d.}
Let $A_1,A_2,\ldots\in G$ be i.i.d. random matrices whose distributions are invariant under both left- and right-multiplication by $K$. Let $\l(k)=\SN(A_1\cdots A_k)$. 
\begin{enumerate} 
\item (Strong law of large numbers) Suppose that the expectation $\E\langle\Cor(A_1),\rho\rangle<\infty$ exists. Then we have almost surely convergence of random vectors in $V$:
$$\frac{\l(k)}{k}\stackrel{a.s.}{\rightarrow}\E\Cor(A_1),\quad k\rightarrow\infty;$$

\item (Central limit theorem) Suppose that the expectation $\E\langle\Cor(A_1),\rho\rangle^2<\infty$ exists. Then we have weak convergence to the multivariate normal distribution
$$\left(\frac{\langle\l(k),\alpha_i\rangle-k\E\langle\Cor(A_1),\alpha_i\rangle}{\sqrt{k}}\right)_{1\le i\le n}\Rightarrow N(0,\Sigma)$$
as $k$ goes to infinity, where $\Sigma=\Cov_{1\le i,j\le n}(\langle\Cor(A),\alpha_i\rangle,\langle\Cor(A),\alpha_j\rangle)$ is the covariance matrix.
\end{enumerate}
\end{thm}

There is a body of work on random walks on Lie structures or root lattices, for instance Montroll \cite{montroll1964random}, Jafarizadeh and Sufiani \cite{jafarizadeh2007investigation}, Grabiner and Magyar \cite{grabiner1993random}, Grabiner \cite{grabiner1999brownian}, Bougerol and Jeulin \cite{bougerol2002paths}, Kabluchko \cite{kabluchko2017convex}. The results of
\Cref{cor: i.i.d.} confirm the idea in \cite[Section 4]{shen2024gaussian} that the connection between singular numbers of products and the sum of corners could be found in all root systems and probability measures. As far as we know, this is the first Gaussian universality result for matrix products that applies to every split reductive group. In particular, for the root system of type $A$, the asymptotic above has the following form: 

\begin{example}
To illustrate \Cref{cor: i.i.d.} concretely, consider the case $G=\SL_{n+1}(F)$. Denote $\l(k)=(\l_1(k),\ldots,\l_{n+1}(k))$. Then we have the following: 
\begin{enumerate} 
\item (Strong law of large numbers) Suppose that for all $1\le i\le n+1$, the expectation $\E|\SN(A_1^{(i)})|<\infty$ exists. Then we have almost surely convergence when $k$ goes to infinity:
$$(\frac{\l_1(k)}{k},\ldots,\frac{\l_{n+1}(k)}{k})\stackrel{a.s.}{\rightarrow}(-\E|\SN(A_1^{(2)})|,\ldots,\\
\E|\SN(A_1^{(n)})|-\E|\SN(A_1^{(n+1)})|,\E|\SN(A_1^{(n+1)})|).$$
Here follow from \Cref{rmk: the name corner}, for all $1\le i\le n+1$, $A_1^{(i)}\in\Mat_{(n-i+2)\times(n+1)}(F)$ denote the corner of $A_1$ consisting of the last $n-i+2$ rows and all $n+1$ columns.

\item (Central limit theorem) Suppose that for all $1\le i\le n+1$, the expectation $\E|\SN(A_1^{(i)})|^2<\infty$ exists. Then we have weak convergence
$$\frac{(\l_1(k)-k\E|\SN(A_1^{(1)})|+k\E|\SN(A_1^{(2)})|,\ldots,\l_{n+1}(k)-k\E|\SN(A_1^{(n+1)})|)}{\sqrt{k}}\Rightarrow N(0,L\Sigma L^T)$$
as $k$ goes to infinity, where $L^T\in\Mat_{(n+1)\times (n+1)}(F)$ is the transpose of the matrix
$$L=\begin{pmatrix}1 & -1 & & &\\ & 1 & -1 & &\\ & & \ddots & \ddots &\\ & & & 1 & -1\\ & & & & 1\end{pmatrix}$$
and $\Sigma=\Cov_{1\le i,j\le n+1}(|\SN(A_1^{(i)})|,|\SN(A_1^{(j)})|)$ is the covariance matrix.
\end{enumerate}
\end{example}

The above example is exactly the Gaussian universality proved in \cite[Theorem 1.1]{shen2024gaussian}; similarly, the case $G=\Sp_{2n}(F)$ coincides with the Gaussian universality proved in \cite[Theorem 4.2]{shen2024gaussian}. The previous proofs for these two examples are purely linear algebra, where $\Cor(\cdot)$ is regarded as the matrix corners obtained by eliminating the rows. While these results rely heavily on explicit matrix structures, this paper realizes $\Cor(\cdot)$ from the perspective of representation theory, which opens pathways for studying random walks on exceptional groups (e.g. $G_2$ or $E_8$) where traditional linear algebraic methods falter. Our asymptotic results in \Cref{cor: i.i.d.} are enabled by a key probabilistic estimate that bounds the discrepancy between singular numbers and corners. This estimate, stated in \Cref{thm: close distance between lambda and nu} below, shows that the difference $\l(k)-\nu(k)$ is negligible, allowing us to equate their asymptotic behaviors:

\begin{thm}\label{thm: close distance between lambda and nu}
Let $\varepsilon,\delta>0$ be fixed, and $A_1,A_2,\ldots\in G$ be independent random matrices. Suppose for all $k\ge 1$, the distribution of $A_k$ is invariant under both left- and right-multiplication by $K$, and the inequality 
$$\mathbf{P}(\SN(A_k)\ne 0)>\delta$$
holds for every $k\ge 1$. Then almost surely, except for finitely many $k$, the following inequality holds:
$$0\le\langle\l(k)-\nu(k),\rho\rangle\le\langle\l(k),\rho\rangle^\varepsilon.$$
\end{thm}

The result in \Cref{thm: close distance between lambda and nu} is stronger and does not require that the matrices $A_1,A_2,\ldots$ share the same distribution. It is worth mentioning that the sequence $\nu(1),\nu(2),\ldots$ can be interpreted as a process involving independent random increments. Specifically, given $A_1,\ldots,A_{k-1}$, we study the random increments $\nu(k)-\nu(k-1)=\Cor(A_1\cdots A_{k-1}A_k)-\Cor(A_1\cdots A_{k-1})$. Since $\Cor$ is right $K$-invariant, and the probability measure of $A_k\in G$ is invariant under left- and right-multiplication of $K$, there is no loss in assuming that $A_1\cdots A_{k-1}\in N\L$. In this case, we have $\Cor(A_1\cdots A_{k-1}A_k)-\Cor(A_1\cdots A_{k-1})=\Cor(A_k)$, which is independent of the previous elements. Therefore, we derive the following equality of distributions (while the two sequences themselves are not identical, they share the same distribution):
\begin{equation}\label{eq: cor product equals sum by distribution}
(\Cor(A_1),\Cor(A_1A_2),\ldots)\stackrel{d}{=}(\Cor(A_1),\Cor(A_1)+\Cor(A_2),\ldots)
\end{equation}
The above fact, together with \Cref{thm: close distance between lambda and nu} allows us to study the asymptotics of the singular numbers of random products by simply adding the independent corners, which gives the proof of \Cref{cor: i.i.d.}.

To draw an analogy of these asymptotic behaviors with previous progress, we closely examine related works on products of random matrices. A prevalent approach in this field is to normalize matrices to a standard form and subsequently analyze their distributions and asymptotics. Inspired by the asymptotic behavior of sums of random variables, such as the law of large numbers and the central limit theorem, Bellman \cite{bellman1954limit}, followed by Furstenberg and Kesten \cite{furstenberg1960products}, examined the limits and fluctuations of singular values in the product sequence $A_1,A_1A_2,\ldots,$ where $A_1,A_2,\ldots$
are random matrices in real or complex fields. One may also view the contributions of Avila-Eskin-Viana \cite{avila2023continuity}, Newman \cite{newman1986distribution}, Vanneste \cite{vanneste2010estimating}, Forrester \cite{forrester2015asymptotics},
Lima-Rahibe \cite{lima1994exact}, Akemann-Ipsen \cite{akemann2015recent}, Liu-Wang-Wang \cite{liu2023lyapunov}. These works often focus on \emph{Lyapunov exponents} from various perspectives, including mathematics and statistical physics. For random complex matrices $A_1,A_2,\ldots$, the $i^{th}$ Lyapunov exponent is defined as
$$\lim_{k\rightarrow\infty}\frac{1}{k}\log(\text{$i^{th}$ largest singular number of $A_1A_2\cdots A_k$})$$
Following the suggestions of \cite{van2021limits} and \cite{shen2024gaussian}, for our split reductive group setting, it is reasonable to regard the $i^{th}$ Lyapunov exponent related to $\alpha_i\in\Delta^+$ for random matrices $A_1,A_2,\ldots\in G$ as the limit
$$\lim_{k\rightarrow\infty}\frac{\langle\l(k),\alpha_i\rangle}{k}=\lim_{k\rightarrow\infty}\frac{\langle\SN(A_1\cdots A_k),\alpha_i\rangle}{k}$$
which emerge from coweight dynamics governed by group-theoretic decompositions (However, in our setting, these Lyapunov exponents no longer decrease). Thus, the strong law of large numbers in \Cref{cor: i.i.d.} gives the Lyapunov exponent $\E\langle\Cor(A_1),\alpha_i\rangle$ related to $\alpha_i\in\Delta^+$ for i.i.d matrices over $G$. 

On the other hand, various techniques have been developed to analyze the non-archimedean case, for which the settings may have little difference. Chhaibi\cite{chhaibi2017non} also discusses the Cartan and Iwasawa decomposition of random matrix products over split reductive groups. Assuming the group has at least one minuscule character, he derives the strong law of large numbers over the Borel subgroup $B$. Our result, however, connects the Cartan and Iwasawa decomposition and therefore gives Gaussian universality, which goes beyond the method in \cite{chhaibi2017non}. Under the setting of Van Peski \cite{van2023p}, the frequency of multiplying new random matrices follows the Poisson process. Brofferio-Schapira \cite{brofferio2011poisson} apply the Oseledets' multiplicative ergodic theorem to prove a law of large numbers for products of i.i.d. random matrices over $\GL_n(\Q_p)$. Finally, a substantial body of work explores random walks on Bruhat-Tits buildings, which provides insights into random walks on $p$-adic groups. See, for example, Cartwright-Woess \cite{cartwright2004isotropic}, Schapira \cite{schapira2009random}, and, in particular, the survey by Parkinson \cite{parkinson2017buildings}.

\subsection{Non-asymptotic results and Satake isomorphism}\label{subsec: Nonasymtotic results and Satake isomorphism}

The asymptotic results in \Cref{subsec: Asymtotic results} are based on substantial work in the realm of non-asymptotic computations, which will be the main concerns of this subsection. For instance, one may study the following questions:

\begin{enumerate}
\item [(Q1)]Let $A,B\in G$ be random. What can we say about the distribution of $\SN(AB)$?

\item [(Q2)]How do singular numbers determine corners? In other words, given the singular number of a random $A\in G$, what can we say about $\Cor(A)$?
\end{enumerate}

We will answer the above questions by applying the Satake isomorphism, which links the spherical Hecke algebra $L(G,K)$ and polynomials $\C[R^\vee]^W$ that are invariant under the Weyl group $W$. Consider the following operations over the Hall-Littlewood polynomials $P_\l(t)$, which form a basis of $\C[R^\vee]^W$: 

\begin{enumerate}
\item (Product Convolution) Since the $P_\l(t)$ form a basis of $\C[R^\vee]^W$, 
$$P_\mu(t) P_\nu(t)=\sum_{\l\in R_+^\vee} c_{\mu,\nu}^\l(t)P_\l(t),\quad \l\in R_+^\vee$$
for some structure coefficients $c_{\mu,\nu}^\l(t)\in\Z[t]$ called the \emph{Littlewood-Richardson coefficients}. 
\item (Polynomial expansion) We expand the polynomial $P_\l(t)$ in a finite linear sum of the form $e^\nu$:
$$P_\l(t)=\sum_{\nu\in R^\vee} u_{\l,\nu}(t)e^{\nu}$$
where $u_{\l,\nu}(t)\in\Z[t]$. 
\end{enumerate}

The polynomials $P_\l(t)$ are, in fact, images of the canonical basis $c_\l$ of $L(G,K)$ under the Satake isomorphism; see \Cref{subsec: Hall-Littlewood polynomials and Weyl groups} and \Cref{subset: Hecke algebra and Satake isomorphism} for further discussion around these polynomials and their properties. We now state the main structural results that cover all root systems and reveal that the singular numbers of products and
corners come from the above operations. The results are, in fact, the generalized version of \cite[Theorem 1.3]{van2021limits}, which concerns the $G=\GL_n(F)$ case. 

\begin{thm}\label{thm: Product process}
(Product process) Let $A,B\in G$ are random with fixed singular numbers $\SN(A)=\mu,\SN(B)=\nu$, invariant under left-and right-multiplication of $K$, we have for any $\l\in R_+^\vee$,
$$\mathbf{P}_{\mu,\nu}^\l=\mathbf{P}(\SN(AB)=\l\mid\SN(A)=\mu,\SN(B)=\nu)=\frac{c_{\mu,\nu}^\l(q^{-1})P_\l(\theta;q^{-1})}{P_\mu(\theta;q^{-1})P_\nu(\theta;q^{-1})}$$
where $\theta$ is the specialization such that $e^{\l}\mapsto q^{\langle\l,\rho\rangle}$, and $q=|\mathfrak{o}/\mathfrak{p}|$ is the order of the residue field. 
\end{thm}

Since singular numbers are invariant under left- and right-multiplication of $K$, it is clear that \Cref{thm: Product process} still holds when one of the matrices is fixed, i.e., let $A\in G$ are random with fixed singular numbers $\SN(A)=\l$, invariant under left- and right-multiplication of $K$, the distribution of singular numbers of $\pi_\l A$ is given by $\mathbf{P}(\SN(\pi_\mu A)=\mu+\nu)=\mathbf{P}_{\mu,\l}^{\mu+\nu}=\frac{c_{\mu,\l}^{\mu+\nu}(q^{-1})P_{\mu+\nu}(\theta;q^{-1})}{P_\mu(\theta;q^{-1})P_\l(\theta;q^{-1})}$. By taking the limit $\mu$ becomes ``very dominant'' and applying \Cref{prop: cor and sn} again, we get the distribution of $\Cor(A)$ as the following:

\begin{thm}\label{thm: Corners process}
(Corners Process) Let $A\in G$ be random with fixed singular numbers $\SN(A)=\l$, invariant under left- and right-multiplication of $K$. Then $\Cor(A)$
is distributed with respect to the law
\begin{equation}\label{eq: conditional distribution of corners}
\mathbf{P}(\nu|\l)=\mathbf{P}(\Cor(A)=\nu)=u_{\l,\nu}(q^{-1})q^{\langle\nu,\rho\rangle}/P_\l(\theta;q^{-1}).
\end{equation}
\end{thm}

As shown in \cite[4.4.4(i)]{bruhat1972groupes}, we always have $\SN(A)\ge\Cor(A)$, hence $\langle\l-\nu,\rho\rangle\ge 0$; conversely, the form of the distribution in \eqref{eq: conditional distribution of corners} indicates that $\langle\l-\nu,\rho\rangle$ tends to be small. Therefore, the corners tend to be not too far away from the singular numbers. Enlightened by this observation, in \Cref{sec: Proof of asymptotic results}, we will apply the Borel-Cantelli lemma to give the rigorous proof of \Cref{thm: close distance between lambda and nu}. Explicit formulas for the probabilities in \Cref{thm: Product process} and \Cref{thm: Corners process} without the appearance of Hall-Littlewood polynomials can be obtained using
the principal specialization formulas for Hall-Littlewood polynomials, see \eqref{Principal specialization formula}. 

\subsection{Plan of paper.} In \Cref{sec: Preliminaries}, we review the split reductive group and probabilistic constructions, representation theory via Satake isomorphism; in \Cref{sec: Proof of nonasymptotic results}, we prove the nonasymptotic results stated in \Cref{subsec: Nonasymtotic results and Satake isomorphism} based on the explicit form provided by Satake isomorphism; in \Cref{sec: Proof of asymptotic results}, we apply the results we get in \Cref{sec: Proof of nonasymptotic results} to prove the asymptotics results in \Cref{subsec: Asymtotic results}. We conclude with open questions beyond left and right $K$-invariant measures.

\section{Preliminaries}\label{sec: Preliminaries}

Let us start with a few paragraphs of backgrounds of matrices over non-archimedean local fields, which are quoted from \cite{shen2024non} and then generalized to arbitrary split reductive groups. Our background in Hall-Littlewood polynomials comes from \cite{macdonald1968spherical}, and connections to the split reductive groups are condensed versions of the results in \cite{macdonald1990orthogonal}. More details of discussions about the typical type $A_n$ can be found in \cite[Section 2, Section 5]{shen2024gaussian}. Also, see \cite[Chapter 1,2]{humphreys1992reflection}, \cite[Chapter 6,7]{sepanski2007compact}, and \cite{humphreys2012linear} for backgrounds in reflection groups, Lie groups, and non-archimedean split reductive groups, respectively.

Fix a non-archimedean local field $F$ with $\mathfrak{o}$ its ring of integers, $\pi$ a generator of the maximal ideal $\mathfrak{p}$, and $q=|\mathfrak{o}/\mathfrak{p}|<\infty$ the order of the residue field. Any nonzero element $x\in F^\times$ could be written as $x=\pi^my$ with $m\in\Z$ and $y\in\mathfrak{o}^\times$. Define $|\cdot|: F \to \R_{\ge 0}$ by setting $|x| = q^{-m}$ for $x$ as before, and $|0|=0$. 

$F$ is equipped with a left- and right-invariant (additive) Haar measure, which is unique if we require the maximal compact subgroup $\mathfrak{o}$ to have measure $1$. The restriction of this measure to $\mathfrak{o}$ is the unique Haar probability measure over $\mathfrak{o}$. As an explicit characterization, its pushforward under any map $r_n:\mathfrak{o}\to\mathfrak{o}/\mathfrak{p}^n$ is the uniform probability measure. It is often useful to view elements of $\mathfrak{o}$ as `power series in $\pi$' $a_0 + a_1 \pi + a_2 \pi^2 + \cdots$, with $a_i\in\{b_0,\ldots,b_{q-1}\}$, a set of representatives of the residue field $\mathfrak{o}/\mathfrak{p}$. Clearly, these specify a coherent sequence of elements of $\mathfrak{o}/\mathfrak{p}^n$ for each $n$. The Haar probability measure then has the alternate explicit description that the $a_i,i\ge 0$ are i.i.d. uniformly picked from $\{b_0,\ldots,b_{q-1}\}$. In addition, $F$ is isomorphic to the ring of the Laurent series in $\pi$, given by the natural isomorphism between the two rings.

\subsection{Split reductive group and random setting}
A \emph{split reductive group} $G=G(F)$ is a reductive group that contains a split maximal torus $T$. The \emph{Weyl group} $W$ of a reductive group $G$ is the quotient group of the normalizer of the maximal torus $T$ by the torus, i.e., $W=N_G(T)/T$. The Weyl group is, in fact, a finite group generated by reflections over roots. 

Similarly to the additive Haar measure on the local field $F$, $G=G(F)$ has a unique left- and right-invariant measure for which the total mass of the maximal compact subgroup $K=G(\mathfrak{o})$ is $1$. The restriction of this measure to $K$, pushes forward to $G(\mathfrak{o}/\mathfrak{p}^n)$ and is the uniform measure on these finite groups. 

The following proposition comes from Cartan decomposition and Iwasawa decomposition, which gives the orbits in $G$ under the multiplication of $K\times K$ and $N\times K$ respectively.

\begin{prop}\label{prop: Cartan and Iwasawa}
\begin{enumerate}
\item (Cartan decomposition) For any $A \in G=G(F)$, there exist $U, V \in K=G(\mathfrak{o})$ such that $A = U\pi_\l V$ where $\l\in R_+^\vee$, and $\pi_\l\in\L$ is the image of $\pi$ under $\l$, viewed as the homomorphism in $\Hom(\G_m,T)$.
\item (Iwasawa decomposition) For any $A\in G$, there exist $U\in N, V \in K$ such that $A = U\pi_\nu V$ where $\nu\in R^\vee$.
\end{enumerate}
\end{prop}

\begin{example}\label{ex: Cartan and Iwasawa for Sp}
The case $G=\SL_{n+1}(F)$ is discussed in \Cref{ex: Cartan and Iwasawa for SL}. Here, we study the group $$G=\Sp_{2n}(F)=\{A\in\GL_{2n}(F)\mid A\begin{pmatrix} & I_n \\ -I_n & \end{pmatrix}A^T=\begin{pmatrix}  & I_n \\ -I_n & \end{pmatrix}\}$$
which is the case of type $C$. Then the Weyl group $W=S_n\rtimes\{1,-1\}^n$ acts naturally on $V=\R^n$. Let $\epsilon_i,1\le i\le n$ denote the standard basis of $\R^n$, then we have
$$\Pi^+=-\Pi^-=\{\epsilon_i-\epsilon_j\mid 1\le i<j\le n\}\cup\{\epsilon_i+\epsilon_j\mid 1\le i\le j\le n\},\rho=(n,n-1,\ldots,1)$$
$$\Delta^+=\{\epsilon_i-\epsilon_{i+1}\mid1\le i\le n-1\}\cup 2\epsilon_n,R=\{(\l_1,\ldots,\l_n)\in\Z^{n}\mid\l_1+\cdots+\l_n \text{ is even}\}, R^\vee=\Z^n$$
$$R^\vee_0=\{(\l_1,\ldots,\l_n)\in\Z^{n}\mid\l_1+\cdots+\l_i\ge 0,\quad\forall 1\le i\le n\}$$
$$R^\vee_+=\{(\l_1,\ldots,\l_n)\in\Z^{n}\mid\l_1\ge\ldots\ge\l_n\ge 0\}.$$
On the other hand, $T$ is the subset of $G$ consisting of diagonal matrices, $\L$ is the subset of $T$ where every entry on the diagonal is some power of $\pi$, 
$$B=\{\begin{pmatrix}B_1 & M\\ 0 & B_1^{-T}\end{pmatrix}\mid B_1\text{ upper triangular}, B_1M^T=MB_1^T\}$$
$$N=\{\begin{pmatrix}N_1 & M\\ 0 & N_1^{-T}\end{pmatrix}\mid N_1\text{ upper triangular with diagonal elements }1, N_1M^T=MN_1^T\}$$
By the Cartan decomposition, there exists $U,V\in\Sp_{2n}(\mathfrak{o})$ such that $A=U\diag_{(2n)\times(2n)}(\pi^{\l_1},\ldots,\pi^{\l_n},\pi^{-\l_1},\ldots,\pi^{-\l_n})V$, where $(\l_1,\ldots,\l_n)\in R_+^\vee$. Also, by the Iwasawa decomposition, there exists $U\in N,V\in\Sp_{2n}(\mathfrak{o})$ such that $A=U\diag_{(2n)\times(2n)}(\pi^{\nu_1},\ldots,\pi^{\nu_n},\pi^{-\nu_1},\ldots,\pi^{-\nu_n})V$, where $(\nu_1,\ldots,\nu_n)\in R^\vee$. Notice that the non-increasing relation $\nu_1\ge\ldots\ge\nu_n\ge 0$ might not hold for the corners.
\end{example}

The restriction of the Haar measure on $G$ to the double coset $K\pi_\l K$ normalized to be a distribution, is the unique $K\times K$-invariant distribution on the subset of $G$ with singular numbers $\lambda$. These distributions are equivalently written as $U \pi_\l V$ where $U,V$ are independently distributed by the Haar probability measure on $K$. More generally, if $U, V \in K$ are Haar distributed and $\l \in R_+^\vee$, then $U \pi_\l V$ is invariant under $K$ acting on both sides, which gives the unique bi-invariant measure with singular numbers given by $\l$. As a corollary, we have the following:

\begin{cor}\label{cor: invariant measure}
Suppose $A\in G$ is random, nonsingular with distribution invariant under the action of $K$ on the left and right. Then $A$ has the form
$$A=U\pi_\l V$$
where the three random matrices on the right-hand side are independent, $U, V\in K$ are Haar distributed, and $\SN(A)=\l$ is randomly distributed over the set $R_+^\vee$.
\end{cor}

\subsection{Hall-Littlewood polynomials and Weyl groups}\label{subsec: Hall-Littlewood polynomials and Weyl groups}

\begin{defi}
Let $\C[R^\vee]$ denote the group algebra over $\C$ of the free abelian group $R^\vee$. For each $\l\in R^\vee$, let $e^\l$ denote the corresponding element of $A$, so that $e^\l e^\mu=e^{\l+\mu}$, $(e^\l)^{-1}=e^{-\l}$, and $e^{0}=1$, the identity element of $\C[R^\vee]$. The group $W$ acts on $R^\vee$ and hence also on $\C[R^\vee]$: $w(e^\l)=e^{w\l}$ for all $w\in W$ and $\l\in R^\vee$. Let $\C[R^\vee]^W$ denote the subalgebra of $\C[R^\vee]$, consisting of elements that are invariant under $W$.
\end{defi}

It is reasonable to view the elements in $\C[R^\vee]$ as formal Laurent series, see the remark on page 10 of \cite{macdonald1990orthogonal}, or \eqref{eq: HL for the Sp case} as an example. Therefore, from now on, we will use the term ``polynomial'' to refer to elements in $\C[R^\vee]$ or $\C[t][R^\vee]$, where $t$ is a parameter.

\begin{defi}
For all $\l\in R_+^\vee$, we denote 
$$W(t)=\sum_{w\in W}t^{n(w)},\quad W_\l(t)=\sum_{w\in W,w\l=\l}t^{n(w)}$$
where $n(w)$ denote the number of positive roots such that $w(\alpha)<0$. It is worth mentioning that $W_\l(t)$ divides $W(t)$, see \cite[Page 21]{humphreys1992reflection}.
\end{defi}

Since each $W$-orbit in $R^\vee$ meets $R_+^\vee$ in exactly one point, it follows that the “monomial symmetric functions”
$$m_\l=\sum_{W/\{w\in W\mid w\l=\l\}}e^{w\l},\quad\l\in R_+^\vee$$
form a $\C$-basis of $\C[R^\vee]^W$. Another basis is provided by the Weyl characters: let 
$$\delta=\prod_{\alpha\in\Pi^+}(e^{\alpha^\vee/2}-e^{-\alpha^\vee/2})=e^{\rho}\prod_{\alpha\in\Pi^+}(1-e^{-\alpha^\vee}).$$
Then $w\delta=\varepsilon(w)\delta$ for all $w\in W$, where $\varepsilon(w)=\det(w)=\pm 1$. For every $\l\in R^\vee$ let
$$\chi_\l=\delta^{-1}\sum_{w\in W}\varepsilon(w)e^{w(\l+\rho)}$$
Then $\chi_\l\in\C[R^\vee]^W$, and the $\chi_\l$ with $\l\in R_+^\vee$ form a $\C$-basis of $\C[R^\vee]^W$. Moreover, we have
$$\chi_\l=m_\l+\sum_{\mu<\l}K_{\l,\mu}m_\mu$$
where $K_{\l,\mu}\in\N$ is the \emph{Kostka number}. If $\l\notin R_+^\vee$, then either $\chi_\l=0$ or else there exists $\mu\in R_+^\vee$ such that $\mu+\rho=w(\l+\rho)$, and in this case $\chi_\l=\varepsilon(w)\chi_\mu$. 

Apart from the two bases $m_\l$ and $\chi_\l$ we discussed above, the definition below provides a basis for $\C[t][R^\vee]^W$, which turns out to be useful throughout our paper.

\begin{defi} Let $W$ be the Weyl group, and $\l\in R_+^\vee$. The \emph{Hall-Littlewood polynomial} $P_\l(t)$ is given by
\begin{equation}\label{eq: HL polynomial}
P_\l(t)=\frac{1}{W_\l(t)}\sum_{w\in W}e^{w\l}\prod_{\alpha\in\Pi^+}\frac{1-te^{-w\alpha^\vee}}{1-e^{-w\alpha^\vee}}=\sum_{w\in W/\{w\in W\mid w\l=\l\}}e^{w\l}\prod_{\alpha\in\Pi^+}\frac{1-te^{-w\alpha^\vee}}{1-e^{-w\alpha^\vee}}.
\end{equation}
\end{defi}
\begin{prop}
The Hall-Littlewood polynomials $P_\l(t)$ satisfies the following properties:
\begin{enumerate}
\item They are monic and have the form
\begin{equation}\label{eq: monomial coefficient}
P_\l(t)=m_\l+\sum_{\mu\in R_+^\vee,\mu<\l} u_{\l,\mu}(t)m_\mu
\end{equation}
here $u_{\l,\mu}(t)\in\Z[t]$.
\item When $\l$ ranges over all the elements in $R_+^\vee$, the Hall-Littlewood polynomials $P_\l(t)$ form a $\C[t]$-basis of $\C[t][R^\vee]^W$.
\end{enumerate}
\end{prop}

In the typical case $W=S_{n+1}$, the notation we use here is slightly different from the previous work in \cite{shen2024gaussian}, where the Hall-Littlewood Laurent polynomials of type $A_n$ are written in the form $P_\l(x_1,\ldots,x_{n+1};t)\in\C[t][x_1^{\pm{1}},\ldots,x_{n+1}^{\pm{1}}]^{S_{n+1}}$, and $\l=(\l_1,\ldots,\l_{n+1})$ satisfies $(\l_1\ge...\ge\l_{n+1})\in\Z^{n+1},\l_1+\cdots+\l_{n+1}=0$. Since we have a canonical homomorphism given by the mapping $\C[R^\vee]^{S_{n+1}}\rightarrow\C[x_1^{\pm{1}},\ldots,x_{n+1}^{\pm{1}}]^{S_{n+1}}, e^\l\mapsto x_1^{\l_1}\cdots x_{n+1}^{\l_{n+1}}$, these two notation coincide. From now on, we sometimes use the form $x_1,\ldots,x_n$ for general root systems and will not distinguish between these two notations. On the other hand, when $W=S_n\rtimes\{1,-1\}^n$ is the Weyl group of $\Sp_{2n}(F)$, for every $\l=(\l_1,\ldots,\l_n)\in R_+^\vee$, we have
\begin{equation}\label{eq: HL for the Sp case}
P_\l(x_1,\ldots,x_n;t)=\frac{1}{W_\l(t)}\sum_{w\in W} x_1^{\l_1}\cdots x_n^{\l_n}\prod_{1\le i<j\le n}\frac{x_i-tx_j}{x_i-x_j}\prod_{1\le i\le n}\frac{x_i-t}{x_i-1}.
\end{equation}
where $W$ is the combination of permutations and inverses $x_i\mapsto x_i^{-1}$ for all $1\le i\le n$.

Another way to express $P_\l(t)$ is to write it as linear sum of Weyl characters: 
\begin{equation}\label{eq: HL as sum of Weyl}
P_\l(t)=\frac{1}{W_\l(t)}\sum_{S\subseteq \Pi^+}(-t)^{|S|}\sum_{w\in W}e^{w(\l-\sum_{\alpha\in S}\alpha)}\prod_{\alpha\in\Pi^+}\frac{1}{1-e^{-w\alpha^\vee}}=\frac{1}{W_\l(t)}\sum_{S\subseteq \Pi^+}(-t)^{|S|}\chi_{\l-\sum_{\alpha\in S}\alpha}
\end{equation}
In particular, when $t=0$, the Hall-Littlewood polynomials $P_\l(t)$ degenerate to the Weyl characters $\chi_\l$.

Since the Hall-Littlewood polynomials $P_\l(t),\l\in R_+^\vee$ form a $\C[t]$-basis of $\C[t][R^\vee]^W$, we have
\begin{equation}\label{eq: LR coef}
P_\mu(t) P_\nu(t)=\sum_{\l\in R_+^\vee}c_{\mu,\nu}^\l(t)P_\l(t)
\end{equation}
for some structure coefficients $c_{\mu,\nu}^\l(t)$, which we call \emph{Littlewood-Richardson coefficient}. 

\begin{prop}\label{prop: LR coef}
The Littlewood-Richardson coefficient $c_{\mu,\nu}^\l(t)$ satisfies the following properties:
\begin{enumerate}
\item $c_{\mu,\nu}^\l(t)\in\Z[t]$. This must be true since the polynomials $P_\l$ are monic and have coefficients in $\Z[t]$;

\item $c_{\mu,\nu}^\l(t)=0$ unless $\mu+\nu\ge\l$. This comes from the fact that the product $P_\mu(t) P_\nu(t)$ has the form
$$P_\mu(t) P_\nu(t)=m_{\mu+\nu}(t)+\text{lower terms}.$$
\end{enumerate}
\end{prop}

One may also expand the form $P_\l(t)$ in a finite linear sum of $e^{\nu},\nu\in R_+^\vee$. It is well known that for all $\nu\in R_+^\vee$ and $w\in W$, we must have $\nu\ge w\nu$, see \cite[Proposition 18, Chapter 6.1]{bourbaki1989lie}. Therefore, the expansion of $P_\l(t)$ is ranged over $\nu\in R^\vee$ such that $\nu\le\l$:
$$P_\l(t)=\sum_{\nu\in R^\vee,\nu\le\l} u_{\l,\nu}(t)e^\nu,\quad\nu\in R^\vee$$
where the $u$ coefficient comes from \eqref{eq: monomial coefficient}. Here, we additionally set $u_{\l,w\nu}(t)=u_{\l,\nu}(t)$ for all $w\in W$ and $\nu\in R_+^\vee$.

\begin{defi}
A \emph{specialization} $\theta$ is an algebra homomorphism from $\C[R^\vee]$ to $\C$. We denote the application of $\theta$ to $f\in\C[R^\vee]$ as $f(\theta)$. In particular, for any $x\in V$, by extending the mappings $e^\l\mapsto q^{\langle\l,x\rangle}$ linearly we obtain a specialization. 
\end{defi}

It is worth mentioning that unlike the typical definition of specialization given in \cite{borodin2013macdonaldprocesses}, the specialization in our paper is defined over $\C[R^\vee]$ instead of the subring of elements invariant under the action of the Weyl group. 

\begin{prop}
(Weyl dimension formula) Let $\theta_0$ be the specialzation $e^\l\mapsto 1$ for all $\l\in R^\vee$. Then we have for all $\l\in R_+^\vee$,
\begin{equation}\label{eq: Weyl dimension formula}
\chi_\l(\theta_0)=\prod_{\alpha\in\Pi^+}\frac{\langle\l+\rho,\alpha^\vee\rangle}{\langle \rho,\alpha^\vee\rangle}
\end{equation}
\end{prop}

\begin{prop}
(Principal specialization formula) Let $\rho=\frac{1}{2}\sum_{\alpha\in\Pi^+}\alpha$ be the Weyl vector, and $\theta$ be the specialization $e^\l\mapsto t^{-\langle\l,\rho\rangle}$ for all $\l\in R^\vee$. Then we have
\begin{equation}\label{Principal specialization formula}
P_\l(\theta;t)=\frac{1}{W_\l(t)}\sum_{w\in W}t^{-\langle w\l,\rho\rangle}\prod_{\alpha\in\Pi^+}\frac{1-t^{1+\langle w\alpha^\vee,\rho\rangle}}{1-t^{\langle w\alpha^\vee,\rho\rangle}}=\frac{W(t)}{W_\l(t)}t^{-\langle\l,\rho\rangle},\quad\forall\l\in R_+^\vee
\end{equation}
\end{prop}

\begin{proof}
\cite[Exercise 6.34]{sepanski2007compact} shows that $\langle \alpha^\vee,\rho\rangle=1$ for all $\alpha\in\Delta^+$; also, for all $w\ne\id$, there exists $\alpha\in\Pi^+$ such that $-w\alpha\in\Delta^+$, and therefore $1-t^{1+\langle w\alpha^\vee,\rho\rangle}=0$. Hence the only nontrivial term in the sum is the case $w=\id$, and
\begin{equation}
P_\l(\theta;t)=\frac{1}{W_\l(t)}t^{-\langle \l,\rho\rangle}\prod_{\alpha\in\Pi^+}\frac{1-t^{1+\langle \alpha^\vee,\rho\rangle}}{1-t^{\langle \alpha^\vee,\rho\rangle}}=\frac{1}{W_\l(t)}t^{-\langle \l,\rho\rangle}\sum_{w\in W}t^{n(w)}=\frac{W(t)}{W_\l(t)}t^{-\langle\l,\rho\rangle}
\end{equation}
where the second equality comes from \cite[Page 84-85]{humphreys1992reflection}.
\end{proof}

\subsection{Hecke algebra and Satake isomorphism}\label{subset: Hecke algebra and Satake isomorphism} Throughout this subsection, $F$ is a local field with $\mathfrak{o}$ its ring of integers, $\pi$ a generator of the maximal ideal $\mathfrak{p}$, and $q=|\mathfrak{o}/\mathfrak{p}|<\infty$ the order of the residue field. Let $\G$ be a split reductive group with Weyl group $W$, $G=G(F),K=G(\mathfrak{o})$ be the maximal compact subgroup.

\begin{defi}\label{defi: Hecke algebra}
Let $L(G,K)$ denote the space of all complex-valued continuous functions of compact support on $G$, which are left- and right-invariant with respect to $K$, i.e., such that $f(k_1xk_2)=f(x)$ for all $x\in G$ and $k_1,k_2\in K$. We define a multiplication on $L(G,K)$ as follows: for $f,g\in L(G,K)$,
$$(f*g)(x)=\int_Gf(xy^{-1})g(y)dy$$
where $dx$ is the unique Haar measure on $G$ such that $K$ has measure $1$. This product is associative and commutative. Also, let $H(G,K)$ denote the subspace of $L(G,K)$ consisting of functions with integer values, which we call the \emph{Hecke algebra} of $G$. 
\end{defi}
\begin{prop}
Every $f\in L(G,K)$ could be written as a finite linear combination of the form $c_\mu$, where $\mu\in R^\vee$, and $c_\mu$ is the characteristic function of the double coset
$$K\pi_\mu K.$$
The similar statement for $H(G,K)$ translates mutatis mutandis.
\end{prop}

\begin{defi}
Let $\mu,\nu\in R_+^\vee$. Then, we define the structure coefficient $G_{\mu\nu}^\l(\mathfrak{o})$ for the expansion of the product $c_\mu*c_\nu$:
\begin{equation}\label{eq: structure coefficient g}
c_\mu*c_\nu=\sum_\l G_{\mu,\nu}^\l(\mathfrak{o})c_\l
\end{equation}
\end{defi}

\begin{thm}
(Satake isomorphism) The mappings
\begin{align}\label{eq: Satake isomorphism}
\begin{split}
\psi: L(G,K)&\rightarrow\C[R^\vee]^W\\
c_\l&\mapsto q^{\langle\lambda,\rho\rangle}P_\l(q^{-1})
\end{split}
\end{align}
extended linearly gives an isomorphism from $L(G,K)$ to $\C[R^\vee]^W$.
\end{thm}

\begin{prop}
$G_{\mu,\nu}^\l(\mathfrak{o})$ has the following properties:
\begin{enumerate}
\item $G_{\mu,\nu}^\l(\mathfrak{o})\in\Z_{\ge 0}$ is a non-negative integer;

\item $G_{\mu,\nu}^\l(\mathfrak{o})$ is a function in $q$, i.e, there exists a polynomial $g_{\mu,\nu}^\l(t)\in\Z[t]$ with degree $\langle\mu+\nu-\l,\rho\rangle$, independent of $\mathfrak{o}$, such that $G_{\mu,\nu}^\l(\mathfrak{o})=g_{\mu,\nu}^\l(q)$.
\end{enumerate}
\end{prop}
\begin{proof}
For the first assertion, notice that
$$G_{\mu,\nu}^{\lambda}(\mathfrak{o})=(c_\mu*c_\nu)(\pi_\lambda)=\int_G c_\mu(y)c_\nu(y^{-1}\pi_\lambda)dy$$
Since $c_\mu(y)$ vanishes for $y$ outside $K\pi_\mu K$, the integration is over this orbit, which we shall write as a disjoint union of left cosets, say
\begin{equation}\label{eq: disjoint left coset}K\pi_\mu K=\bigsqcup_j y_jK\quad (y_j\in K\pi_\mu)\end{equation}
Therefore, we have
$$G_{\mu,\nu}^{\lambda}(\mathfrak{o})=\sum_j\int_{y_jK}c_\nu(y^{-1}\pi_\lambda)dy=\sum_jc_\nu(y_j^{-1}\pi_\lambda)$$
since $K$ has measure $1$. Hence $G_{\mu,\nu}^{\lambda}(\mathfrak{o})$ is equal to the number of $j$ such that $y_j^{-1}\pi_\lambda\in K\pi_\nu K$, which is a non-negative integer; Also, \eqref{eq: LR coef}, \eqref{eq: structure coefficient g}, \eqref{eq: Satake isomorphism} together connect the coefficient $G_{\mu,\nu}^\l(\mathfrak{o})$ which corresponds to the Hecke algebra, and $c_{\mu,\nu}^\l(q^{-1})$ which corresponds to the group algebra $\C[R^\vee]^W$:
\begin{equation}\label{eq: G and c}
G_{\mu,\nu}^\l(\mathfrak{o})=q^{\langle\mu+\nu-\l,\rho\rangle}c_{\mu,\nu}^\l(q^{-1})
\end{equation}
where right hand side is purely over $q$, independent of $\mathfrak{o}$. Since this form is an integer for every $q$ as a power of prime, the degree of the polynomial $c_{\mu,\nu}^\l(q^{-1})\in\Z[q^{-1}]$ must be less than $\langle\mu+\nu-\l,\rho\rangle$, and therefore the explicit form
in \eqref{eq: G and c} is a polynomial in $q$ with integer coefficients. 
\end{proof}

From now on, we always write $g_{\mu,\nu}^\l(q)$ for the structure coefficient instead of $G_{\mu,\nu}^\l(\mathfrak{o})$.

\section{Proof of non-asymptotic results}\label{sec: Proof of nonasymptotic results}

In this section, we prove \Cref{thm: Product process} and \Cref{thm: Corners process} based on the background provided in the Preliminaries section. We will postpone our proof of the asymptotic results to the next section after we finish studying \Cref{thm: Product process} and \Cref{thm: Corners process}. Let us start with the following lemma.

\begin{lemma}\label{lem: transition step}
For all $\lambda\in R_+^\vee$, let $$V(K\pi_\lambda K)=\int_G c_\lambda(x)dx$$ 
denote the volume of the orbit $K\pi_\l K$ (i.e., the number of left cosets in \eqref{eq: disjoint left coset}), where $dx$ is the $G$-invariant measure on $G$ normalized by $\int_K dx=1$. Suppose $A,B\in G$ are random with fixed singular numbers $\SN(A)=\mu,\SN(B)=\nu$, invariant under left-and right multiplication of $K$. Then the probability $\mathbf{P}_{\mu,\nu}^{\lambda}:=\mathbf{P}(\SN(AB)=\lambda)$ from \Cref{thm: Product process} has the form
$$\mathbf{P}_{\mu,\nu}^{\lambda}=\frac{g_{\mu,\nu}^{\lambda}(q)V(K\pi_\lambda K)}{V(K\pi_\mu K) V(K\pi_\nu K)}.$$
\end{lemma}

\begin{proof}[Proof of \Cref{lem: transition step}]
Consider the integral 
$$\mathcal{I}=\int_{G\times G}c_\nu(x)c_\mu(y)c_\lambda(xy)dxdy.$$
On the one hand, we have
\begin{equation}\label{eq: x and y integral}\mathcal{I}=V(K\pi_\nu K)\int_G c_\mu(y)c_\lambda(\pi_\nu y)dy=\mathbf{P}_{\mu,\nu}^{\lambda}V(K\pi_\mu K) V(K\pi_\nu K).\end{equation}
The second equality holds because the set on which $c_\mu(y)=1$ has measure $V(K\pi_\mu K)$ (since each coset has measure $1$), and the proportion of this set on which $c_\lambda(\pi_\nu y)=1$ equals $\mathbf{P}_{\mu,\nu}^{\lambda}$. On the other hand, set $z=xy$. Since the measure is $G$-invariant, we have
\begin{align}\label{eq: z and y integral}
\begin{split}
\mathcal{I}&=\int_{G\times G}c_\lambda(z)c_\mu(y)c_\nu(zy^{-1})dzdy\\
&=V(K\pi_\lambda K)\int_G c_\mu(y)c_\nu(\pi_\lambda y^{-1})dzdy=g_{\mu,\nu}^{\lambda}(q)V(K\pi_\lambda K)
\end{split}
\end{align}
where the last equality is by definition of $g_{\mu,\nu}^{\lambda}(q)$. The two results from \eqref{eq: x and y integral} and \eqref{eq: z and y integral} together give the proof.
\end{proof}

With the above preparation, we turn back to our proof of \Cref{thm: Product process} and \Cref{thm: Corners process}.

\begin{proof}[Proof of \Cref{thm: Product process}]
For the product process, by \cite[(4.5')(4.7)]{macdonald1968spherical}, we have for all $\l\in  R_+^\vee$,
\begin{equation}\label{eq: volumn of double coset}
V(K\pi_\lambda K)=q^{2\langle\lambda,\rho\rangle}\frac{W(t)}{W_\l(t)}
\end{equation}
Apply \Cref{lem: transition step} and the results from \eqref{eq: G and c} and \eqref{eq: volumn of double coset}, we have
\begin{align}
\begin{split}
\mathbf{P}_{\mu,\nu}^{\lambda}&=\frac{g_{\mu,\nu}^{\lambda}(q)V(K\pi_\lambda K)}{V(K\pi_\mu K) V(K\pi_\nu K)}\\
&=q^{\langle\l-\mu-\nu,\rho\rangle}c_{\mu,\nu}^\l(q^{-1})\frac{W_\mu(q^{-1})W_\nu(q^{-1})}{W_\l(q^{-1})W(q^{-1})}\\
&=\frac{c_{\mu,\nu}^\l(q^{-1})P_\l(\theta;q^{-1})}{P_\mu(\theta;q^{-1})P_\nu(\theta;q^{-1})}
\end{split}
\end{align}
where the last row comes from the principal specialization formula in \eqref{Principal specialization formula}. 

In the end, the property of the coefficient $c_{\mu,\nu}^\l(q^{-1})$ in \Cref{prop: LR coef} implies that $\mathbf{P}_{\mu,\nu}^\l=0$ unless $\mu+\nu\ge\l$, or equivalently $\SN(A)+\SN(B)\ge\SN(AB)$ almost surely.
\end{proof}

\begin{rmk}\label{rmk: bound of SN of product}
It is well known that for any $A,B\in G$, we always have $\SN(A)+\SN(B)\ge\SN(AB)$, see \cite[4.4.4 (iii)]{bruhat1972groupes}. Combining our result with \Cref{prop: LR coef} verifies this property, since $c_{\mu,\nu}^\l(t)=0$ unless $\mu+\nu\ge\l$.
\end{rmk}

\begin{example}
Consider the example $G=\SL_{n+1}(F)$. Denote 
$$\mu=(\mu_1,\ldots,\mu_{n+1}),\nu=(\nu_1,\ldots,\nu_{n+1})$$ 
where $(\mu_1\ge\ldots\ge\mu_{n+1}),(\nu_1\ge\ldots\ge\nu_{n+1})\in\Z^{n+1}$, and $\mu_1+\cdots+\mu_{n+1}=\nu_1+\cdots+\nu_{n+1}=0$. We write $$P_\mu(x_1,\ldots,x_{n+1};q^{-1})P_\nu(x_1,\ldots,x_{n+1};q^{-1})=\sum_{\l}c_{\mu,\nu}^\l(q^{-1}) P_\l(x_1,\ldots,x_{n+1};q^{-1})$$ 
as the linear sum of Hall-Littlewood polynomials. In this case, $\theta$ is the specialization sending $x_i$ to $q^{\frac{n-2i+2}{2}}$ for all $1\le i\le n+1$. Therefore, we have
\begin{align*}
\mathbf{P}_{\mu,\nu}^{\l}&=\frac{c_{\mu,\nu}^\l(q^{-1})P_\l(q^{n/2},\ldots,q^{-n/2};q^{-1})}{P_\mu(q^{n/2},\ldots,q^{-n/2};q^{-1})P_\nu(q^{n/2},\ldots,q^{-n/2};q^{-1})}\\
&=\frac{c_{\mu,\nu}^\l(q^{-1})P_\l(1,\ldots,q^{-n};q^{-1})}{P_\mu(1,\ldots,q^{-n};q^{-1})P_\nu(1,\ldots,q^{-n};q^{-1})}
\end{align*}
where the second rows holds because $P_\mu(x_1,\ldots,x_{n+1};t),P_\nu(x_1,\ldots,x_{n+1};t),P_\l(x_1,\ldots,x_{n+1};t)$ are homogenous of degree zero, see \cite[Lemma 2.2]{van2021limits}. Hence, the form we get coincides with the previous result from \cite[Theorem 1.3]{van2021limits}.
\end{example}

\begin{example}
Consider the case $G=\Sp_{2n}(F)$. Denote $\mu=(\mu_1,\ldots,\mu_n),\nu=(\nu_1,\ldots,\nu_n)\in R_+^\vee$. We write $$P_\mu(x_1,\ldots,x_n;q^{-1})P_\nu(x_1,\ldots,x_n;q^{-1})=\sum_{\l}c_{\mu,\nu}^\l(q^{-1}) P_\l(x_1,\ldots,x_n;q^{-1})$$ 
as the linear sum of Hall-Littlewood polynomials of the form in \eqref{eq: HL for the Sp case}. In this case, $\theta$ is the specialization sending $x_i$ to $q^{n-i+1}$ for all $1\le i\le n$. Therefore, we have
\begin{equation*}
\mathbf{P}_{\mu,\nu}^{\l}=\frac{c_{\mu,\nu}^\l(q^{-1})P_\l(q^n,\ldots,q;q^{-1})}{P_\mu(q^n,\ldots,q;q^{-1})P_\nu(q^n,\ldots,q;q^{-1})}
\end{equation*}
\end{example}

The following proposition shows that given any fixed $A\in G$, the corners $\Cor(A)$ could be obtained by observing the singular numbers $\SN(\pi_\mu A)$ where $\mu$ is ``sufficiently dominant''. 

\begin{prop}\label{prop: cor and sn}
Let $A\in G$ be fixed. Then there exists $M>0$ such that for all $\mu\in R_+^\vee$ that satisfies $$\langle\mu,\alpha\rangle>M,\quad\forall\alpha\in\Pi^+$$
we must have $\Cor(A)=\SN(\pi_\mu A)-\mu$.
\end{prop}

\begin{proof}
Since $\Cor:G\rightarrow R^\vee$ is invariant under the right multiplication of $K$, there is no loss to assume $A\in N\L$. In this case, we have $A\pi_{\Cor(A)}^{-1}=\pi_0+\sum_{\alpha\in\Pi^+}c_\alpha\mathfrak{g}_{\alpha}\in N$, where $c_\alpha\in F$ denotes the coefficient of $\mathfrak{g}_\alpha$. The adjoint representation over the weight spaces $\mathfrak{g}_{\alpha}$ takes the form
$$\pi_\mu \mathfrak{g}_{\alpha}\pi_\mu^{-1}=\pi^{\langle\mu,\alpha\rangle}\mathfrak{g}_{\alpha},\quad\forall \mu\in R^\vee,\alpha\in\Phi$$
We pick $M>0$ such that $\pi^Mc_\alpha\in\mathfrak{o}$ for all $\alpha\in\Pi^+$. Then for all $\mu\in R_+^\vee$ that satisfies $\langle\mu,\alpha\rangle>M,\forall\alpha\in\Pi^+$, we must have 
$$\pi_\mu A\pi_{\Cor(A)}^{-1}\pi_\mu^{-1}=\pi_0+\sum_{\alpha\in\Pi^+}\pi^{\langle\mu,\alpha\rangle}c_\alpha\mathfrak{g}_{\alpha}\in N(\mathfrak{o})\subset K,\quad\pi_\mu A\in K\pi_{\mu+\Cor(A)}.$$
Therefore, we have $\SN(\pi_\mu A)=\mu+\Cor(A)$, which gives the proof.
\end{proof}

\begin{rmk}\label{rmk: SN bigger than Cor}
It is well known that for any $A\in G$, we always have $\SN(A)\ge\Cor(A)$, see \cite[4.4.4 (i)]{bruhat1972groupes}. Combining our result with \Cref{rmk: bound of SN of product} verifies this property, since $\Cor(A)=\SN(\pi_\mu A)-\mu\le\SN(A)$.
\end{rmk}

\begin{proof}[Proof of \Cref{thm: Corners process}]

Consider the decomposition of the double coset $K\pi_\l K$ as the finite joint union of the left cosets, say $K\pi_\l K=\bigsqcup_j y_jK$, where $y_j\in N\L$ for all $1\le j\le V(K\pi_\l K)$. Then $\mathbf{P}(\nu\mid\l)$ becomes the proportion of those $y_j$ such that $\Cor(y_j)=\nu$. Since there are only finitely many $j$, \Cref{prop: cor and sn} implies that there exists a uniform $M>0$ such that for all $\mu\in R_+^\vee$ that satisfies $\langle\mu,\alpha\rangle>M,\forall\alpha\in\Pi^+$, we must have $\Cor(y_j)=\SN(\pi_\mu y_j)-\mu$, i.e., $\mathbf{P}(\nu\mid\l)=\mathbf{P}_{\l,\mu}^{\mu+\nu}$. We take a coweight $\mu\in R_+^\vee$ that meets this requirement and also satisfies
$$\langle\mu+\nu,\alpha_i\rangle>0$$ for all $u_{\l,\nu}(t)\ne 0$ and $1\le i\le n$. Then we apply the form in \eqref{eq: HL polynomial}:
\begin{align}
\begin{split}
P_\l(t)P_\mu(t)&=\sum_{w\in W}w(e^{\mu}P_\l(t)\prod_{\alpha\in\Pi^+}\frac{1-te^{-\alpha^\vee}}{1-e^{-\alpha^\vee}})\\
&=\sum_\nu u_{\l,\nu}(t)\sum_{w\in W}w(e^{\mu+\nu}\prod_{\alpha\in\Pi^+}\frac{1-te^{-\alpha^\vee}}{1-e^{-\alpha^\vee}})\\
&=\sum_{\nu}u_{\l,\nu}(t)P_{\mu+\nu}(t).
\end{split}
\end{align}
The third row is because the only $w\in W$ that satisfies $w(\mu+\nu)=\mu+\nu$ is $w=\id$. Hence $c_{\l,\mu}^{\mu+\nu}(t)=u_{\l,\nu}(t)$, and
\begin{equation}
\mathbf{P}(\nu\mid\l)=\mathbf{P}_{\l,\mu}^{\mu+\nu}=\frac{u_{\l,\nu}(q^{-1})P_{\mu+\nu}(\theta;q^{-1})}{P_\l(\theta;q^{-1})P_\mu(\theta;q^{-1})}=u_{\l,\nu}(q^{-1})q^{\langle\nu,\rho\rangle}/P_\l(\theta;q^{-1})
\end{equation}
where the last equality comes from the principal specialization formula \eqref{Principal specialization formula}.
\end{proof}

\section{Proof of asymptotic results}\label{sec: Proof of asymptotic results}

In this section, we will prove \Cref{thm: close distance between lambda and nu}, and based on this estimate of discrepancy, we show the Gaussian universality in \Cref{cor: i.i.d.}. We begin with the following technical lemmas.

\begin{lemma}\label{lem: lower and upper bound of the corners}
Let $\l\in R_+^\vee$. Then for all $\nu\in R^\vee$ such that $\mathbf{P}(\nu\mid\l)>0$, we have
$$-\langle\l,\rho\rangle\le\langle \nu,\rho\rangle\le\langle\l,\rho\rangle$$
\end{lemma}

\begin{proof} 
Let $w\in W$ such that $v=w\nu_+$, where $\nu_+\in R_+^\vee$. Then we also have $\mathbf{P}(\nu_+\mid\l)>0$. On one hand, by \cite[Proposition 18, Chapter 6.1]{bourbaki1989lie}, we always have $\l\ge\nu_+\ge \nu$, thus $\langle \nu,\rho\rangle\le\langle\l,\rho\rangle$; On the other hand, there exists $w_0\in W$ that $w_0\Pi^+=\Pi^-,w_0\rho=-\rho$. Therefore, we have
$$\langle \l+\nu,\rho\rangle=\langle \l-w_0\nu,\rho\rangle\ge\langle \nu_+-w_0\nu,\rho\rangle\ge 0$$
which ends the proof.
\end{proof}

\begin{lemma}\label{lem: exponential decrease}
\begin{enumerate}
\item For all $\l\in R_+^\vee,\nu\in R^\vee$, we have $u_{\l,\nu}(q^{-1})\ge 0$. Also, we have $\sum_{\nu\in R^\vee}u_{\l,\nu}(q^{-1})=O(\langle\l,\rho\rangle^{|\Pi|})$ when $\langle\l,\rho\rangle\rightarrow\infty$. Here by the $O$ symbol we mean that there exists a fixed constant $C>0$ such that $\sum_{\nu\in R^\vee}u_{\l,\nu}(q^{-1})\le C\langle\l,\rho\rangle^{|\Pi|}$ for all $\l\in R_+^\vee$;

\item For every $C,\varepsilon>0$, as $\langle\l,\rho\rangle\rightarrow\infty$, we have
$$\sum_{\nu\in R^\vee,\langle\l-\nu,\rho\rangle\ge C\langle\l,\rho\rangle^{\varepsilon}}\mathbf{P}(\nu\mid\l)= O(\langle\l,\rho\rangle^{|\Pi|}q^{-C\langle\l,\rho\rangle^{\varepsilon}}).$$   
\end{enumerate}
\end{lemma}

\begin{proof}
On the one hand, the coefficient $u_{\l,\nu}(q^{-1})$ must be non-negative because the probability $\mathbf{P}(\nu\mid\l)\ge 0$ is nonnegative; On the other hand, let $\theta_0$ be the specialization $e^\nu\mapsto 1$ for all $\nu\in R^\vee$, we have
\begin{align*}
\sum_{\nu\in R^\vee}u_{\l,\nu}(q^{-1})&=\sum_{\nu\in R^\vee}u_{\l,\nu}(q^{-1})e_\nu(\theta_0)\\
&=P_\l(\theta_0;q^{-1})\\
&=\frac{1}{W_\l(q^{-1})}\sum_{S\subseteq \Pi^+}(-q)^{-|S|}\chi_{\l-\sum_{\alpha\in S}\alpha}(\theta_0)=O(\langle\l,\rho\rangle^{|\Pi|}).
\end{align*}
where the last row comes from \eqref{eq: HL as sum of Weyl} and Weyl dimension formula \eqref{eq: Weyl dimension formula}.

In the end, based on the above estimate, we have
\begin{align*}\sum_{\nu\in R^\vee,\langle\l-\nu,\rho\rangle\ge C\langle\l,\rho\rangle^{\varepsilon}}\mathbf{P}(\nu\mid\l)&=\sum_{\nu\in R^\vee,\langle\l-\nu,\rho\rangle\ge C\langle\l,\rho\rangle^{\varepsilon}}\frac{u_{\l,\nu}(q^{-1})q^{\langle\nu,\rho\rangle}}{P_\l(\theta;q^{-1})}\\
&\le\sum_{\nu\in R^\vee,\langle\l-\nu,\rho\rangle\ge C\langle\l,\rho\rangle^{\varepsilon}} u_{\l,\nu}(q^{-1})q^{-C\langle\l,\rho\rangle^{\varepsilon}}\\
&=O(\langle\l,\rho\rangle^{|\Pi|}q^{-C\langle\l,\rho\rangle^{\varepsilon}})
\end{align*}
where $\theta$ is the specialization such that $e^\nu\mapsto q^{\langle\nu,\rho\rangle}$, and the second row comes from the principal specialization formula \eqref{Principal specialization formula}. This ends the proof.
\end{proof}

\begin{lemma}\label{lem: lower bound for corner expectation}
There exists $\eta>0$ such that for any $\l\ne 0\in R_+^\vee$, and $A\in G$ be random with $\SN(A)=\l$ fixed and distribution invariant under the left- and right-multiplication of $K$, we have 
$$\E\langle\Cor(A),\rho\rangle=\sum_\nu\mathbf{P}(\nu\mid\l)\langle\nu,\rho\rangle=\sum_\nu \frac{1}{P_\l(\theta;q^{-1})}u_{\l,\nu}(q^{-1})q^{\langle\nu,\rho\rangle}\langle\nu,\rho\rangle>\eta$$
\end{lemma}

\begin{proof}
First, we show $\E\langle\Cor(A),\rho\rangle>0$ for all $\l\ne 0\in R_+^\vee$. Denote $S_\nu=\sum_{w\in W}q^{\langle w\nu,\rho\rangle}\langle w\nu,\rho\rangle=\sum_{w\in W}q^{\langle \nu,w\rho\rangle}\langle \nu,w\rho\rangle$, we have
$$\E\langle\Cor(A),\rho\rangle=\sum_{\nu\in R_+^\vee}\frac{u_{\l,\nu}(q^{-1})}{|\{w\in W\mid w\nu=\nu\}|P_\l(\theta;q^{-1})}S_\nu.$$
The coefficient is non-negative for every $S_\nu$ and must be strictly positive when $\nu=\l$. Hence, we only need to prove $S_\nu>0$ for all nonzero $\nu\in R_+^\vee$. Notice that there exists $w_0\in W$ that $w_0\Pi^+=\Pi^-,w_0\rho=-\rho$. Therefore, we have
\begin{align*}
2S_\nu&=\sum_{w\in W}(q^{\langle \nu,w\rho\rangle}\langle \nu,w\rho\rangle+q^{\langle \nu,ww_0\rho\rangle}\langle \nu,ww_0\rho\rangle)\\
&=\sum_{w\in W}(q^{\langle \nu,w\rho\rangle}-q^{\langle \nu,-w\rho\rangle})\langle \nu,w\rho\rangle>0
\end{align*}
which must be positive. Next, apply the result from \Cref{lem: lower and upper bound of the corners} and \Cref{lem: exponential decrease}, we know that 
$$\sum_\nu\mathbf{P}(\nu\mid\l)\langle\nu,\rho\rangle=(1+o(1))\langle\l,\rho\rangle,\quad\langle\l,\rho\rangle\rightarrow\infty$$
Therefore, we are done.
\end{proof}

\begin{prop}\label{prop: linear increase of cor sum}
Let $A_1,A_2,\ldots\in G$ be the same as in \Cref{thm: close distance between lambda and nu}. Let $\eta$ be the same as in \Cref{lem: lower bound for corner expectation}. Then almost surely, we have
$$\inf_{k\ge 1}\frac{\langle\Cor(A_1)+\cdots+\Cor(A_k),\rho\rangle}{k}\ge\delta\eta.$$
\end{prop}

\begin{proof}
For every $k\ge 1$, set 
$$X_k=\max\{-k^{1/3},\min\{\langle\Cor(A_k),\rho\rangle,k^{1/3}\}\}=\begin{cases}
k^{1/3} & \langle\Cor(A_k),\rho\rangle>k^{1/3}\\
-k^{1/3} & \langle\Cor(A_k),\rho\rangle<-k^{1/3}\\
\langle\Cor(A_k),\rho\rangle & \text{ else}
\end{cases}$$
Then apply \Cref{lem: exponential decrease}, we know
\begin{align*}
\sum_{k\ge 1}\mathbf{P}(X_k>\langle\Cor(A_k),\rho\rangle)&=\sum_{k\ge 1}\mathbf{P}(\langle\Cor(A_k),\rho\rangle<-k^{1/3})\\
&=\sum_{k\ge 1}\mathbf{P}(\langle\SN(A_k),\rho\rangle>k^{1/3},\langle\Cor(A_k),\rho\rangle<-k^{1/3})\\
&\le\sum_{k\ge 1}O(k^{|\Pi|/3}q^{-2k^{1/3}})<\infty.
\end{align*}
where the second row comes from \Cref{lem: lower and upper bound of the corners}, and the last row comes from \Cref{lem: exponential decrease}. Hence, by Borel-Cantelli lemma we know that almost surely, there only exists finitely many $k\ge 1$ such that $X_k>\langle\Cor(A_k),\rho\rangle$. On the other hand, we have
$$\E(X_k)=\sum_{\l\ne 0\in R_+^\vee}\mathbf{P}(\SN(A_k)=\l)\sum_{\nu\in R^\vee}\mathbf{P}(\nu\mid\l)\max\{-k^{1/3},\min\{\langle\nu,\rho\rangle,k^{1/3}\}\}$$
By \Cref{lem: exponential decrease}, we have $\sum_{\nu\in R^\vee}\mathbf{P}(\nu\mid\l)\max\{-k^{1/3},\min\{\langle\nu,\rho\rangle,k^{1/3}\}\}=(1+o(1))k^{1/3}$ as $\langle\l,\rho\rangle$ goes to infinity. Hence for sufficiently large $n$, we have $\E(X_k)\ge\eta\sum_{\l\ne 0\in R_+^\vee}\mathbf{P}(\SN(A_k)=\l)>\delta\eta$. Since $X_k\in[-k^{1/3},k^{1/3}]$, we must have $\sum_{k\ge 1}\frac{\var(X_k)}{k^2}<\infty$, so we can apply Kolmogrov's criterion of the strong law of large numbers, which shows that almost surely,
$$\inf_{k\ge 1}\frac{X_1+\cdots+X_k}{k}\ge\delta\eta.$$
Thus the similar statement also holds when replacing $X_k$ by $\langle\Cor(A_k),\rho\rangle$.
\end{proof}

Now, we may turn back to our proof of \Cref{thm: close distance between lambda and nu}.

\begin{proof}[Proof of \Cref{thm: close distance between lambda and nu}]
For all $k\ge 1$, let $\l(k)=\SN(A_1\cdots A_k)$ and $\nu(k)=\Cor(A_1\cdots A_k)$. As stated in \Cref{rmk: SN bigger than Cor}, we always have $\l(k)\ge\nu(k)$. Denote the event
$$C_k=\{\langle\l(k)-\nu(k),\rho\rangle\ge\langle\l(k),\rho\rangle^{\varepsilon},\frac{\langle\l(k),\rho\rangle}{k}\ge\frac{1}{2}\delta\eta\},\quad k\ge 1.$$

On one hand, by \Cref{lem: exponential decrease}, we have
\begin{equation*}
\sum_{k\ge 1}\mathbf{P}(C_k)=\sum_{k\ge 1}O((\frac{1}{2}k\delta\eta)^{|\Pi^+|} q^{-(\frac{1}{2}k\delta\eta)^{\varepsilon}})<\infty
\end{equation*}
Therefore, the Borel-Cantelli Lemma implies that almost surely, there are only finitely many $ k\ge 1$ in which the event $C_k$ happens. 

On the other hand, as stated in \eqref{eq: cor product equals sum by distribution}, the random sequence $\nu(1),\nu(2),\ldots$ has the same distribution as $\langle\Cor(A_1),\rho\rangle,\langle\Cor(A_1)+\Cor(A_2),\rho\rangle,\ldots$. Also, by  \Cref{rmk: SN bigger than Cor}, we always have $\l(k)\ge\nu(k)$. Hence \Cref{prop: linear increase of cor sum}, almost surely, except for only finitely many $k\ge 1$, we have $\frac{\langle\l(k),\rho\rangle}{k}\ge\frac{\langle\nu(k),\rho\rangle}{k}>\frac{1}{2}\delta\eta$. Together, these two observations give the proof.
\end{proof}

\begin{rmk}
Apart from the proof above, we give another perspective that could help explain why \Cref{thm: close distance between lambda and nu} should be correct, which is a generalized version of the heuristic in \cite[Section 1]{shen2024gaussian}. Suppose $A_1,\ldots,A_k$ are already given, and we write $\l(k)=\SN(A_1\cdots A_k)$. We are interested in determining the distribution of
$\SN(A_1\cdots A_{k+1})$. While this product may seem challenging to handle directly,
we take advantage of the fact that the distribution is invariant under left- and right-multiplication of $K$. This allows us to transfer the matrix $A_1\cdots A_k$ to its normal form under Cartan decomposition and then analyze the singular numbers of the resulting matrix.

Thus, we examine the singular numbers $\l(k+1)=\SN(\pi_{\l(k)}A_{k+1})\in R_+^\vee$. As $k$ goes to infinity, we expect to see $\l(k)$ becomes ``very dominant", i.e., for all positive root $\alpha\in\Pi^+$, $\langle\l(k),\alpha\rangle$ is very large. Then by \Cref{prop: cor and sn}, the increment of singular numbers when multiplying the new matrix $A_{k+1}$ is just $\Cor(A_{k+1})$, i.e., $\l(k+1)$ is approximately equal to $\l(k)+\Cor(A_{k+1})$. Therefore, the asymptotic of singular numbers of matrix products could be directly inherited from the sum of corners. 
\end{rmk}

Based on the asymptotic result in \Cref{thm: close distance between lambda and nu}, the i.i.d. case of \Cref{cor: i.i.d.} becomes a particular circumstance, and the Gaussian universality emerges as the following proof.

\begin{proof}[Proof of \Cref{cor: i.i.d.}]
For all $k\ge 1$, let $\l(k)=\SN(A_1\cdots A_k)$ and $\nu(k)=\Cor(A_1\cdots A_k)$ be the same as above. \eqref{eq: cor product equals sum by distribution} indicates that the sequence $\nu(1),\nu(2),\ldots$ is the sum of i.i.d. variables. Therefore, it has the following limits and fluctuations: 

\begin{enumerate} 
\item (Strong law of large numbers) Suppose that the expectation $\E\langle\Cor(A_1),\rho\rangle<\infty$ exists. Then we have almost surely convergence of random vectors in $V$:
$$\frac{\nu(k)}{k}\stackrel{a.s.}{\rightarrow}\E\Cor(A_1),\quad k\rightarrow\infty;$$

\item (Central limit theorem) Suppose that the expectation $\E\langle\Cor(A_1),\rho\rangle^2<\infty$ exists. Then we have weak convergence to the multivariate normal distribution
$$\left(\frac{\langle\nu(k),\alpha_i\rangle-k\E\langle\Cor(A_1),\alpha_i\rangle}{\sqrt{k}}\right)_{1\le i\le n}\Rightarrow N(0,\Sigma)$$
as $k$ goes to infinity, where $\Sigma=\Cov_{1\le i,j\le n}(\langle\Cor(A_1),\alpha_i\rangle,\langle\Cor(A_1),\alpha_j\rangle)$ is the covariance matrix.
\end{enumerate}

If $\mathbf{P}(\SN(A_1)=0)=1$, there is nothing to prove; Otherwise, suppose that there exists a nonzero $\l\in R_+^\vee$ such that $\mathbf{P}(\SN(A_1)=\l)=\delta>0$. In this case, we can apply the result in \Cref{thm: close distance between lambda and nu}. Observe that for any $\langle\nu,\rho\rangle\ge 1$, the inequality
$\langle\l,\rho\rangle\ge\langle\nu,\rho\rangle+2\langle\nu,\rho\rangle^{1/4}$ implies the inequality $\langle\l,\rho\rangle-\langle\l,\rho\rangle^{1/4}\ge\langle\nu,\rho\rangle$, we have almost surely, the events $$\langle\l(k)-\nu(k),\rho\rangle\ge 2\langle\nu(k),\rho\rangle^{1/4},\quad k\ge 1$$
occur only finitely many times, i.e., except for finitely many $k$, when we express the difference $\l(k)-\nu(k)$ as the non-negative sum of positive coroots, every coefficient is less than $2\langle\nu(k),\rho\rangle^{1/4}$, which is negligible under strong law of large numbers and central limit theorem.
\end{proof}

Here are some possible ways that might lead to further research.

\begin{enumerate}
\item The current version of \Cref{thm: close distance between lambda and nu} requires $\mathbf{P}(\SN(A_k)\ne 0)=\mathbf{P}(A_k\in K)>\delta$ for all $k\ge 1$. However, when we multiply a new matrix $A_k$ to the sequence such that $\SN(A_k)=0$, both singular numbers and corners remain where they are. Therefore, we would expect the same property to hold for weaker conditions. Specifically, we conjecture that as long as
$$\sum_{k\ge 1}\mathbf{P}(\SN(A_k)\ne 0)=\infty$$
then almost surely, except for finitely many $k$, we have $0\le\langle\l(k)-\nu(k),\rho\rangle\le\langle\l(k),\rho\rangle^\varepsilon$. By the (second) Borel-Cantelli lemma, this requirement is equivalent to the assertion that almost surely, there exist infinitely many $k\ge 1$ such that $\SN(A_k)\ne 0$, and therefore cannot be further strengthened. (Otherwise, singular numbers and corners of the sequence of products stop somewhere and never move anymore.)
\item Throughout our paper, we focus on the case that the probability measures over $G$ are left- and right-invariant under the multiplication of $K$, for which the Hecke algebra becomes a valuable tool. So what can we say for more general probability measures when Satake isomorphism no longer works? In such cases, how should we study the stochastic process of matrix products?
\end{enumerate}


\end{document}